\documentclass[journal,transmag]{IEEEtran}

\usepackage{graphicx}
\usepackage{amsmath}
\usepackage{amssymb}

\ifCLASSINFOpdf
\else
\fi
\begin{document}
\title{Voltage Control in the Presence of Photovoltaic Systems }


\author{\IEEEauthorblockN{Ashkan Zeinalzadeh\IEEEauthorrefmark{1},
Reza Ghorbani \IEEEauthorrefmark{2},
James Yee\IEEEauthorrefmark{2}
}
\IEEEauthorblockA{\IEEEauthorrefmark{1}Center for Sustainable Energy at Notre Dame, IN 46556, USA}
\IEEEauthorblockA{\IEEEauthorrefmark{2} College of Engineering, University of Hawai`i at Manoa, Honolulu, HI 96822, USA}
}

%\markboth{Journal of \LaTeX\ Class Files,~Vol.~11, No.~4, December~2012}%
%{Shell \MakeLowercase{\textit{et al.}}: Bare Demo of IEEEtran.cls for Journals}



% for Transactions on Magnetics papers, we must declare the abstract and
% index terms PRIOR to the title within the \IEEEtitleabstractindextext
% IEEEtran command as these need to go into the title area created by
\IEEEtitleabstractindextext{%
\begin{abstract}

In the past decade, the landscape of energy production has had to shift to accommodate renewables. Which, unlike fossil fuels, are subject to frequent fluctuations; potentially destabilizing grid operators. With the continued demand for solar PV system installation, there is a pressing need for utilities to regulate the voltages at the low voltage distribution grids. We develop a stochastic model for voltage rise as a function of injected power into the grid. This model is formed as a linear combination of gamma random variables. This is achieved by finding sparse bases and clustering the data into subsets by its correlation with those bases, and fitting a gamma distribution within each subset. We are concerned with modeling voltage rise, while taking into account sparse events in the voltage. We use sparse singular value decomposition (SVD) with $\ell_1$ penalty to model sparse voltage rise. More randomness and disorder in voltage rises was observed at the point of common coupling (PCC) of the PV systems with greater line impedance. A controller to regulate the voltage at the PCC of a single phase solar PV system is presented. This controller uses the stochastic model to minimize the risk of voltage rise. Simulation results confirm that this voltage controller can lower the voltage more effectively than conventional voltage regulators during periods of high solar PV output.
\end{abstract}

}

\maketitle

\IEEEdisplaynontitleabstractindextext

% \fi
\IEEEpeerreviewmaketitle

\section{Introduction}

\begin{figure}[h!]
\includegraphics[width=8cm]{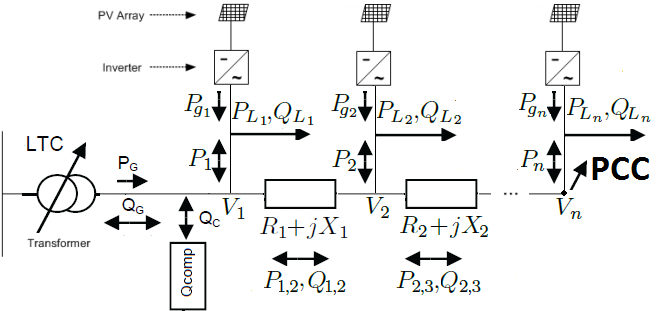}
\caption{Low voltage radial distribution grid.}
\label{fig:distribution}
\end{figure}

\IEEEPARstart This paper is an extended version of work published in \cite{conferenceversion}, which does not provide the design of the voltage regulator. Power injected  to  the  grid  causes  a  voltage  rise  at  the point of common coupling (PCC) due to the inverter and line impedance. In order to determine the maximum capacity of solar photovoltaic (PV) systems on a distribution grid, voltage rise from the penetration of solar energies \cite{c15}-\cite{c17777} must be quantified. A maximum allowed PV capacity is set by the utilities to avoid voltage violations or insufficient load on the grid \cite{c10}. Voltage rise is considered in the design of PV systems, e.g., the inverter disconnects from the grid if the output voltage exceeds its operating limit. Electric grids with high penetration of solar PV systems are highly dependent on solar energies to satisfy the demanded load. The tripping of solar PV systems can lead to periods of low voltage in the distribution grid. The stochastic nature of the loads and renewable energies makes it increasingly difficult to adjust voltage regulators and capacitor banks in the presence of high penetration of renewable energies. A description of a low voltage, radial distribution grid is provided in Figure~1. A point of common coupling (PCC), shown in Figure~\ref{fig:distribution}, is a point on the distribution grid where a consumer is connected.

Let $P_{L_k}$ and $Q_{L_k}$ be the active and reactive load of the $k$th consumer. $P_{g_k}$ is defined as the output of the PV system at the $k$th consumer. $Q_c$ is the reactive compensator connected to the distribution line. $R_k$ and $X_k$ are, respectively, the resistance and reactance of the line between consumers $k$ and $k+1$.
The aggregate power of consumer $k$ is defined as $P_k=P_{L_k}-P_{g_k}$ for all $k=1,...,n$. Negative and positive values for $P_k$ are indicated by the downward and upward arrows for $P_k$ respectively.
Negative power values for $P_k$ imply injection of power into the grid and positive power values represent power flowing in the opposite direction, from the grid to the consumer. Let $P_{k,k+1}$ and $Q_{k,k+1}$ be the active and reactive power flow on the line between consumer $k$ and $k+1$. The voltage at the PCC of the $(k+1)$th consumer can be approximately calculated as:
\begin{equation}\label{eq:volt}
V_{k+1}^2\approx V_k^2-2(R_k P_{k,k+1}+X_k Q_{k,k+1}).
\end{equation}

$P_{k,k+1}$ and $Q_{k,k+1}$ are a function of the load and PV output of consumers $k+1,...,n$. The PCC voltage of a PV system is dependent on both active and reactive power injected into the grid. In this study, PV systems only inject active power into the grid, owing that the amount of reactive load is negligible for residential consumers. Therefore, we do not consider reactive power in this voltage analysis. It is evident from (\ref{eq:volt}) that the voltage at each PCC is affected by neighboring PV system outputs and loads. Therefore, quantifying  PV generation at a single PCC cannot explain the voltage variation at the corresponding PCC. Furthermore, PV system outputs and loads are random variables, and voltage variations can be modeled as random processes.

In this work each consumer can estimate the voltage rise of its respective PCC by simply observing its aggregate power (load minus PV generation). The voltage rise at the PCC of each consumer is modeled only as a function of local information, i.e. using the aggregate power of the consumer. This stochastic model provides an estimation of the voltage rise in the distribution grid without communication among the consumers.

The voltage rise, caused by the PV systems, can negatively effect the operation of voltage regulators, e.g. capacitor banks and transformer load tap changers (LTCs). Specifically, a higher number of tap changes resulting from PV system variability can lead to LTC degradation. In this work, a controller is designed to adjust the LTC using the stochastic model. It is shown in the simulation that the
stochastic model enables us to adjust the LTC more efficiently to decrease the voltage at times of high PV solar output.

In the model for the voltage rise, sparse events are taken into account. Voltage rise has been modeled as a linear combination of Gamma random variables. This model determines the probability of voltage rise as a function of active power at the PCC. Due to variations in the load and PV system output, large increases in PCC voltage can occur within a short time interval. Sparse singular value decomposition (SVD) with $\ell_1$ penalty has been used to model sparse voltage rises. The data related to an instance of high voltage can account for a small portion of the whole dataset; consequently, the direct application of standard SVD to the data would be ineffectual. An algorithm is developed to find the bases for the power data taken from the PCC of consumers using $\ell_1$-singular value decomposition ($\ell_1$-SVD). The data of the consumers on a radial distribution grid is stacked in a single matrix. In the stacked matrix, the rows denote the days, and columns denote the time. We consider the right singular vectors of the $\ell_1$-SVD of the stacked matrix as the bases for the rows. A sparse basis implies that a small number of rows of the stacked matrix, or equivalently a small number of days, are correlated to that basis. The right singular vectors provide a feature for the active power variation of the consumers. Data elements for the days that are correlated with a given basis, are grouped in the corresponding subset for that basis. The voltage data corresponding to each subset of data has been modeled separately as a Gamma random variable. The parameters of the Gamma distributions are evaluated using maximum likelihood estimation (MLE). Finally, the voltage has been modeled as a linear combination of Gamma random variables.

The voltage rise caused by PV systems has been studied in several comparable works. Unlike this paper, these studies do not consider the dynamic effect of the neighbours' loads and PV system outputs on the voltage of a single PCC. A survey study on the voltage rises caused by PV solar outputs is presented below.

The authors of \cite{c7} employ battery storage to avoid voltage violations and alleviate unwanted PV shut-off. It is not cost efficient for individual consumers to own energy storage units \cite{c77}. The authors in \cite{c8}, \cite{Sera} and \cite{ashkanreactive} propose a control algorithm to regulate the voltage of a PV system. They control the voltage by injecting reactive power into the grid. However, lowering the voltage by injecting reactive power into the grid can increase the loss. The authors have not considered grid loss in their analysis. The application of energy storage to mitigate voltage rises was also proposed in \cite{Alam} and \cite{Ueda}. The excess energy from the solar PV systems, is used to charge a distributed energy storage unit during midday, while the stored energy is used to reduce the peak load in the evening. This application of energy storage to regulate voltage may not be economical. The authors in \cite{Ali} study the effect of solar irradiation and penetration levels of PV systems on the voltage quality of a distribution grid in the UK.

The rest of the paper is structured as follows. In Section~\ref{SectionII} the data used in the analysis is described. In Section~\ref{SectionV} an algorithm is described to find appropriate bases to cluster the active power data. These bases enable us to cluster the sparse events into separate subsets. In Section~\ref{SectionVI} a voltage model is developed as a linear combination of Gamma random variables. In Section~\ref{Sectionft}, a comparison of real data and the developed voltage model is presented. In Section~\ref{se:vr} a voltage regulator is developed based on the stochastic model of voltage and the numerical results are presented. Finally, concluding remarks are provided in Section~\ref{SectionVII}.

\section{Data}\label{SectionII}

The data obtained from five consumers, each equipped with a PV system, consists of their respective voltages and average powers.
The consumers are located at different distances from the distribution transformer. The PV system capacities and the line impedance of  consumers are given in Table~I below. The line impedance is defined as the impedance from the PCC of the consumer to the distribution transformer.

\begin{table}[h!]
\caption{Line impedance (OHM) and PV system capacity (kW)}
\label{table:tb0}
\begin{center}
\begin{tabular}{|c|c|c|c|c|c|}
\hline
  Capacity (kW) &1.9&3.9&7.3&11.6&9.2\\
\hline
   Impedance (OHM)  &0.077&0.060&0.053&0.025&0.011\\
%\hline
%   Average daily energy/SC  &4.67&5.13&3.96&5.17&4.98\\
\hline
\end{tabular}
\end{center}
\end{table}

\begin{figure}[h!]
\framebox{\parbox{3in}{
\includegraphics[width=3in]{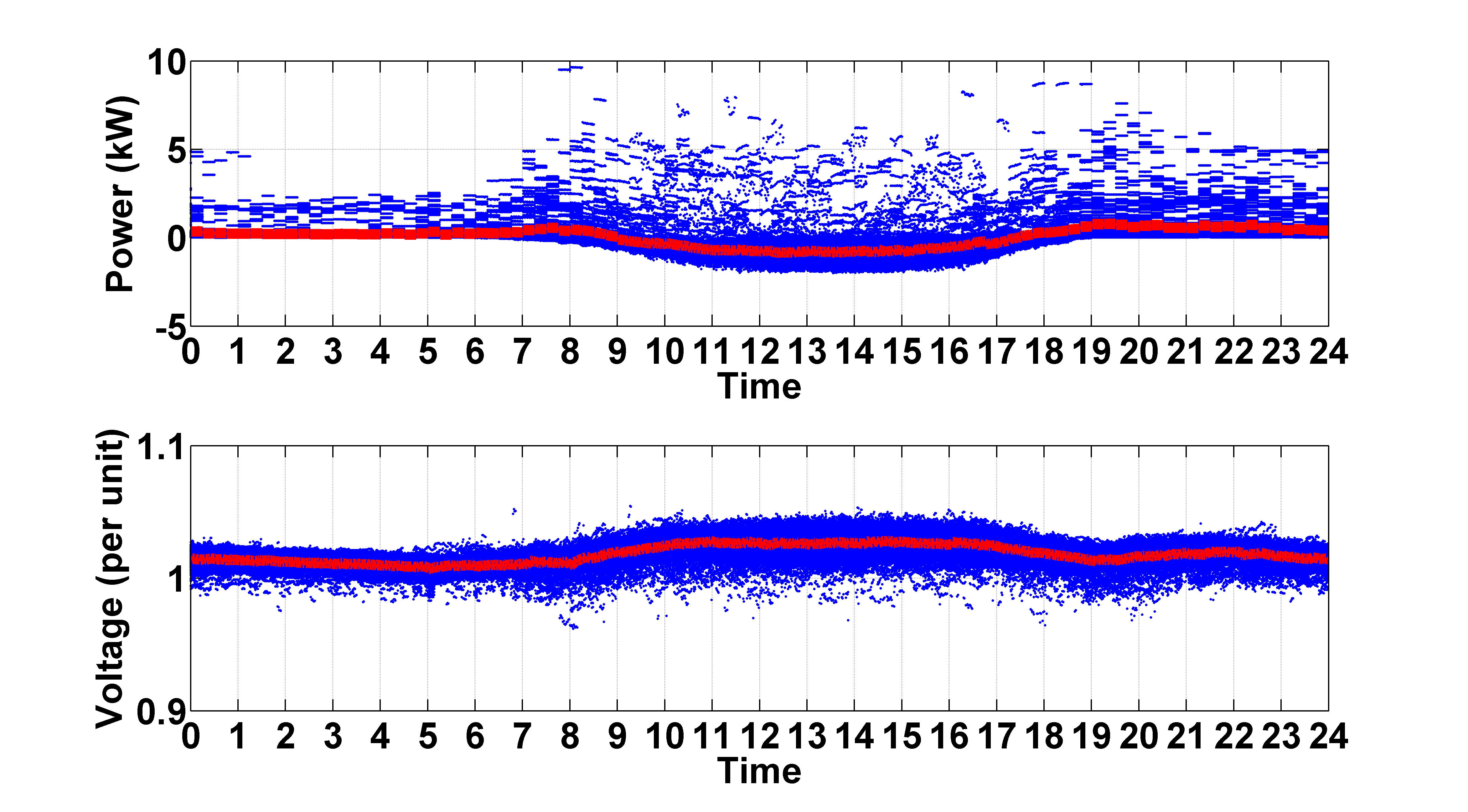}}}
\caption{Power and voltage versus time for the PV system capacity $1.9$ kW.}
\label{fig:povo1}
\end{figure}

\begin{figure}[h!]
\framebox{\parbox{3in}{
\includegraphics[width=3in]{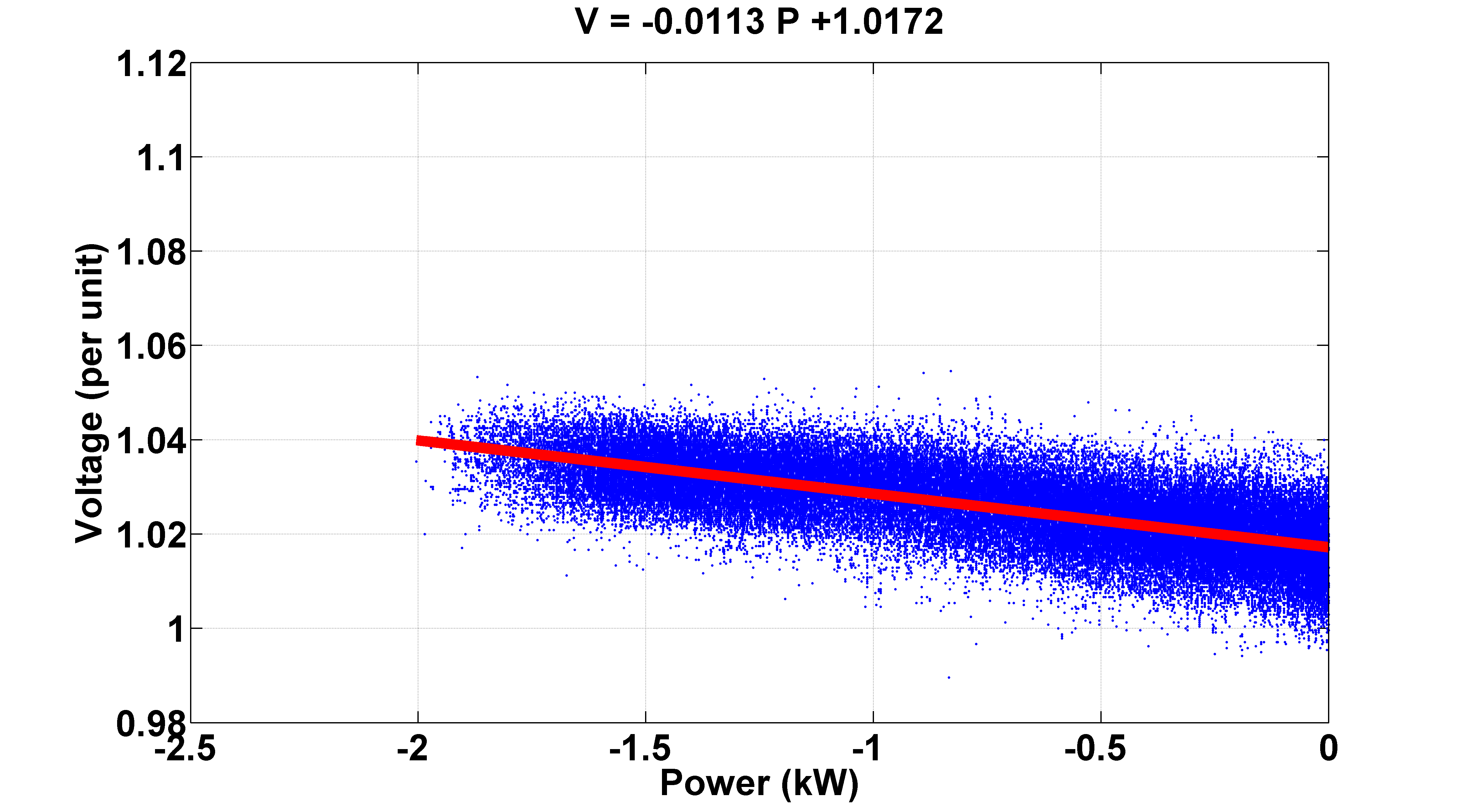}}}
\caption{Voltage versus power for the PV system capacity $1.9$ kW.}
\label{fig:vvp1}
\end{figure}

The PV units in this study are equipped with measurement devices that report the net power (load minus PV generation) and RMS voltage each minute. In Figure~\ref{fig:povo1}, power and voltage at the PCC for the PV system $1.9$ kW are plotted. The blue curves are the averages of the data at each minute for $160$ days. Voltage rise is observed from $10:00$ a.m. to $5:00$ p.m., due to active power injected into the grid from PV systems. In Figure~\ref{fig:vvp1}, the voltage versus the power for the PV system $1.9$ kW is plotted. Each point presents the voltage value for the corresponding injected power. The blue line is the best linear fit to the data using least-squares regression. From Figure~\ref{fig:povo1} and \ref{fig:vvp1}, it is observed that the voltage does not rise consistently when more power is injected into the grid at a single PCC, due to the stochastic nature of the load and PV generation of the neighbors. It is observable that \emph{on average} injecting more power into the grid, causes a voltage rise.

The goal of the next two sections is to develop a stochastic model to represent the voltage at a PCC. The voltage at a PCC is a complex function of the voltages and net powers at the residence as well at neighboring residences. We develop a method to parsimoniously represent the net power data at $5$ residences in Hawaii. The approach will be to modify the method of principal components using singular value decomposition. In addition, a clustering method is used to group or summarize data that is highly correlated (collinear). This is continued in the Section~\ref{SectionVI} to model the voltage at a PCC as a linear combination of the net powers and a set of Gamma random variables.

\section{New Basis for the Power Data} \label{SectionV}

For each residence, the data is stored in the two $m \times n$ matrices $P$ and $V$. The element $p_{ij}$ in row $i$ and column $j$ of the matrix $P$ is the average injected power on day $i$, $i = 1, 2, . . . , m$ and minute $j$, $j = 1, 2, . . . , n$. In the dataset, $m = 160$ days and $n = 24 \times 60$ is the number of minutes in a day. Whenever $p_{ij} > 0$, the residence draws power from the grid. Conversely, when $p_{ij} < 0$, the residence injects power into the grid. The element $v_{ij}$ of the matrix $V$ is the largest RMS voltage in the interval $[j, j + 1)$ on day $i$.

Let $D(t)$ and $G(t)$ be respectively two stochastic processes that represent the load and solar power generated in the interval $[t, t + 1)$. These processes can be modeled as Markov processes. Consider the reshaping of the matrix $P$ into a row vector consisting of the first row, followed by the second row, etc. The resulting row vector would be a sample path of the random process $D(t) - G(t)$. Another stochastic process of interest is the average injected power in the minute $[j, j + 1)$ averaged over $n$ days. For this stochastic process, $\frac{1}{n} \sum_{i=1}^{n} p_{ij}$ is a sample path. The voltage at the PCC is a deterministic function of the injected average powers. Since the average injected powers are random variables, the voltages are also random variables. For $k = 1, 2, . . . , N$, let $P_k$ and $V_k$ respectively be the average injected powers and voltages at residence $k$. In the dataset, $N = 5$. The $m \times n$ matrix $H$ will represent the dataset. $H$ is formed by starting with $P_1$ and appending $P_2$, $P_3$ and so forth. The net power measurements for residence $k$ are stored in the $m \times n$ matrix $P_k$ where the rows correspond to the days and the columns correspond to the minutes of a day. The net power data for the $N$ residences is stored in the matrix $H=\begin{bmatrix} P_1^T & P_2^T & \ldots & P_N^T \end{bmatrix}^T$. The matrix $H$ is formed by appending $P_2$ below $P_1$ then appending $P_3$ below $P_2$ and so forth. The approach will be to capture the essence of the data in $H$ with another matrix with far fewer rows.

Our approach begins by simply applying the singular value decomposition (SVD) theorem. Let $m$ be the number of rows in $H$. This constitutes the total number of days in which data was collected for $N$ houses. By the SVD theorem, $H$ can be expressed as $H = X\Sigma Y^T$ where $X$ is $m \times m$, $Y$ is $n \times n$, $\Sigma$ is $m \times n$, with rank $r$ and positive diagonal entries $\sigma_1 \geq \ldots \geq \sigma_r > 0$. The columns $x_1, \ldots, x_m$ of $X$ are called left singular vectors and $X$ is orthonormal, i.e. $X^TX = I$. The columns $y_1, . . . , y_n$ of $Y$ are categorized as right singular vectors and $Y$ is orthonormal, i.e. $Y^T Y = I$. The $\sigma_i$ are called singular values.

The SVD has the useful property that the splitting of the matrix $H$ into rank-one parts ordered by their singular values is given by
\begin{equation}\label{eq:outprdct}
H = x_1 \sigma_1 y_1^T+\ldots+x_r \sigma_r y_r^T,
\end{equation}
where $r$ is the rank. In the computational work reported in section V, the numerical values of $\sigma_k$ for $k \geq 5$ were very small. Consequently, in a pure principal components approach, one only needs about $4$ terms in the expansion in (\ref{eq:outprdct}). In each execution of the \textbf{Procedure} below, one left and one right singular vector are found. In the initial step (step~$0$), the starting vectors are left and right singular vectors of the SVD given by the $H$ under consideration.

The next part of the \textbf{Procedure} is to modify the left and right singular vectors. Let $t$ denote the iteration count when repeating of steps $1$ and $2$. Let $x^t$ and $y^t$ denote the left and right singular vectors in iteration $t$. In step $1$, suppose the current right singular vector is $y^t$. The left singular vector, $x^{t+1}$ is updated by solving the optimization problem
\begin{align}\label{eq:optimization}
\min_{x} \Big{\{}  \| H-x(y^t)^{T} \|_{F}^2 + \alpha \|x\|_{1} \Big{\}}.
\end{align}

The Frobenius norm of an $m \times n$ matrix $A$ is $\|A\|_{F}= \sqrt{\sum_{i=1}^{m} \sum_{j=1}^{n} a_{ij}^{2}}$. In (\ref{eq:optimization}), $H$ is either the original net power data matrix or the residual data matrix in later iterations. In the first part of (\ref{eq:optimization}), the selection of $x^{t+1}$ is aimed at minimizing the Euclidean distance between $H$ and $x^{t+1}(y_t)^T$. The second part of (\ref{eq:optimization}) is directed towards selecting a $x^{t+1}$ which is sparse. For any vector $x$, the $\ell_0$ norm $\|x\|_{0}$ is equal to the number of nonzero elements of $x$. For example, if $\|x\|_{0}=K$, then $x$ has $K$ nonzero elements and $m-K$ elements are equal to zero. Consequently, as $\|u\|_{0}$ is made small, $u$ becomes sparser. The parameter $\alpha$ is a numerical constant which is selected experimentally. As $\alpha$ is increased, more emphasis is placed on sparsity. In the computational work, it was found that $\alpha = 0.05$ was a good choice.

With the $\ell_0$ norm $\|x\|_{0}$ in the objective function, the problem is combinatorial in nature and commonly formulated as an integer programming problem. However, this approach results in a computationally expensive procedure. In the field of compressive sensing, it is common to replace the $\ell_0$ norm with the $\ell_1$ norm $\|x\|_{1}=\sum_{i} \mid x_i \mid$. Although the $\ell_1$ norm is nonlinear, there is a trick in optimization to formulate this as a linear program. Consequently, (\ref{eq:optimization}) becomes a convex optimization problem for which there are many efficient algorithms.\\
\\
\textbf{Procedure(H)}:
\begin{itemize}
\item \textbf{Step~$0$}: Find $x_1$ and $y_1$, the left and right singular vectors, corresponding to the largest singular value of $H$ by using the Matlab command svd.

\item \textbf{Step~$1$}: For $y^t$, find $x^{t+1}$ by solving the optimization problem
\begin{align}\label{eq:step1}
x^{t+1}=\arg\min_{x} \Big{\{}  \| H-x(y^t)^{T} \|_{F}^2 + \alpha \|x\|_{1} \Big{\}}.
\end{align}
\item \textbf{Step~$2$}: For $x^{t+1}$, calculate $y^{t+1}$ by
\begin{align}\label{eq:step2}
y^{t+1}=\frac{H^T x^{t+1}}{\|x^{t+1}\|}.
\end{align}
If $\| x^{t+1}-x^{t} \| < \epsilon$, then go to the final step~$3$. Otherwise, update the iteration counter and return to step~$1$.
\item \textbf{Step~$3$}: At convergence, normalize the singular vectors and singular value by
\begin{align}\label{eq:step3}
x_1=\frac{x^t}{\|x^t\|}, \,\ y_1=\frac{y^t}{\|y^t\|}, \,\ \sigma_1=x_1^T H y_1.
\end{align}
\end{itemize}

In each application of the \textbf{Procedure}, the row vector $y_1^T$ is found as an additional basis vector for $H$. In the algorithm below, let $H_t$ be the residual matrix at the beginning of iteration $t$. In step~$0$, the residual matrix at the beginning of iteration~$1$ is $H_1=H$. In step~$1$, the rows of $H$ that are highly correlated with $y_1^T$ are found. Let $\mathcal{R}=\{1,2,\ldots,m\}$ denote the set of row indices of $H$. Let $\mathcal{C}_1 \subseteq \mathcal{R}$ be the indices of the rows that are highly correlated with $y_1^T$. The criterion applied is for each $i \in \mathcal{R}$ where the $i$th element of $x_1$ is positive, say $x_{1i} > 0$, then $i$ is included in $\mathcal{C}_1$.

In step~$2$ where $\|x_1\|_0=1$, the basis vector $y_1^T$ is not significant. The rows that have not been included in some $\mathcal{C}_t$ are collected together as the remaining cluster. Let $C$ be the number of clusters or groups. Then $C-1$ is the number of iterations of the algorithm. In the computational results, $C = 3$. Collectively, the set of clusters is exhaustive, i.e. $\mathcal{R}  \subseteq \mathcal{C}_1 \cup \ldots \cup \mathcal{C}_C$, but not necessarily mutually exclusive, i.e. $\mathcal{C}_1 \cap \ldots \cap \mathcal{C}_C \neq  \varnothing$.\\
\\
\textbf{Algorithm:}

\begin{itemize}
\item \textbf{Step~$0$}: Set the iteration counter to $t=1$. Set $H_t = H$.

\item \textbf{Step~$1$}: For $H_t$, find $x_1$, $y_1$, $\sigma_1$ by using \textbf{Procedure($H_t$)}. Using $x_1$ and $y_1$, determine $\mathcal{C}_t$. If $\|x_1\|_0=1$ then go to the final step~$2$. Otherwise, calculate the residual $H_{t+1}=H_t-\sigma_1 x_1 y_1^T$. Increase the iteration counter to $t+1$ and repeat this step.

\item \textbf{Step~$2$}: Include in $\mathcal{C}_{t+1}$ those rows of $H$ that have not been included in $\mathcal{C}_1,\ldots,\mathcal{C}_t$.
\end{itemize}

In the next section, the voltage variation is modeled as a function of injected power into the grid.

\section{Modeling the voltage at a PCC} \label{SectionVI}
In Section~\ref{SectionV}, the data of the days that were correlated with one basis were grouped in one subset.
In this Section, the voltage data corresponding to each subset is modeled separately as a Gamma random variable.
Then, the voltage at each PCC is written as the linear combination of Gamma random variables.
The coefficient of each random variable is the ratio of the number of days in the subset corresponding to the sum of the sizes of all subsets.\\

For each residence, data is stored in two $m \times n$ matrices $P$ and $V$. The elements $p_{ij}$ and $v_{ij}$ in row $i$ and column $j$ of the matrix $P$ and $V$ are the corresponding power and voltage on the same day $i$ and minute $j$. The least square regression parameter ${\beta}$ is given as
\begin{equation}
{\beta}= \Big(\sum_{i}^{m}\sum_{j}^{n} p_{ij}^2\Big)^{-1} \sum_{i=1}^{m} \sum_{j=1}^{n} v_{ij} p_{ij}.
\end{equation}
Assume the power data $P$ is grouped into the different subsets $\mathcal{C}= \{\mathcal{C}_1,...,\mathcal{C}_C \}$ as explained in the Section~\ref{SectionV}. Note that the subsets may not be disjoint, namely because the power data of one day may correlate with multiple bases. Each subset can include the data of multiple days. The sets $I_k^+$ and $I_k^+$ for $k=1,\ldots,C$ are defined as follows,
\begin{equation}
I_k^+=\{ij | i \in \mathcal{C}_k \,\ \text{and} \,\ v_{ij}-\beta p_{ij} \geq 0\},
\end{equation}
\begin{equation}
I_k^-=\{ij | i \in \mathcal{C}_k \,\ \text{and} \,\ v_{ij}-\beta p_{ij} < 0\}.
\end{equation}
Let $n_k^s$ be the size of ${I}_k^{s}$ and $\pi_k^s$ is defined as follows
\begin{equation}
\pi_k^s=\frac{n_k^s}{\sum_{i=1}^{C}(n_i^++n_i^-)}.
\end{equation}
Let $s \in \{+,-\}$, $u_k^s$ is a random variable that is gamma-distributed with shape $\lambda_k^s$ and scale $\theta_k^s$ parameters,
\begin{equation}\label{eq:nu}
u_k^s   \sim \text{Gamma}(\lambda_k^s,\theta_k^s).
\end{equation}
The parameters $\lambda_k^s$ and $\theta_k^s$ are chosen by implementing maximum likelihood estimation such that the random variable $u_k^s$,
represents the best estimate for the subset $\{\frac{|v_{ij}-\beta p_{ij}|}{\pi_k^{s}} \}_{ij \in {I}_k^{s}}$. The log-likelihood function for $n_k^s$ independent and identically distributed observations $\{\frac{|v_{ij}-\beta p_{ij}|}{\pi_k^{s}} \}_{ij \in {I}_k^{s}}$ is

\begin{align}\label{eq:likelihood}
&\ell(\lambda_k^s,\theta_k^s)=(\lambda_k^s-1) \sum_{ij \in {I}_k^{s}} \log(\frac{|v_{ij}-\beta p_{ij}|}{\pi_k^{s}})\\ \notag
&-\sum_{ij \in {I}_k^{s}} \frac{|v_{ij}-\beta p_{ij}|}{\pi_k^{s}}\theta_k^s+n_k^s \lambda_k^s \log(\theta)-n_k^s \log(\Gamma(\lambda_k^s)).
\end{align}

Let $\acute{\Gamma}$ be the derivative of gamma function ${\Gamma}$. (\ref{eq:eq1}) and (\ref{eq:eq2}) are obtained by differentiating the log likelihood function (\ref{eq:likelihood}) with respect to parameters $\lambda_k^s$ and $\theta_k^s$ respectively and setting to zero,
\begin{equation}\label{eq:eq1}
{n}_k^{s} \log \theta_k^{s}+\sum_{ij \in {I}_k^{s}} \log (\frac{|v_{ij}-\beta p_{ij}|}{\pi_k^{s}})-{n}_k^{s} \frac{\acute{\Gamma}(\lambda_k^{s})}{\Gamma(\lambda_k^{s})}=0,
\end{equation}
\begin{equation}\label{eq:eq2}
{\theta_k^{s}}=\frac{{n}_k^{s} \pi_k^{s} \lambda_k^{s}}{ \sum_{ij \in {I}_k^{s}} |v_{ij}-\beta p_{ij}|}.
\end{equation}
(\ref{eq:eq1}) and (\ref{eq:eq2}) are conditions for finding the parameters $\lambda_k^{s}$ and $\theta_k^{s}$ which maximize the log likelihood function (\ref{eq:likelihood}).
Note that (\ref{eq:eq1}) is a nonlinear equation, it can be solved using an iterative numerical procedure such as Newton iteration method or by using the Matlab command gamfit.

Let $p$ be the measured net power at the PCC of the consumer. We are interested in modeling the voltage at the PCC of the consumer as a function of measured net power $p$.
This model includes a linear combination of gamma random variables $\{u_k^+\}_{k=1}^{C}$ and $\{u_k^-\}_{k=1}^{C}$ with coefficients $\{\pi_k^+\}_{k=1}^{C}$ and $\{\pi_k^-\}_{k=1}^{C}$ respectively.
\begin{equation}\label{eq:V}
v=\beta p+\sum_{k=1}^{C} \pi_k^+ u_k^+-\pi_k^- u_k^-,
\end{equation}
\begin{equation}
\sum_{k=1}^{C} \pi_k^++\pi_k^- =1,
\end{equation}
where $\beta$ is the regression parameter.

The Gamma distribution represents a family of shapes. For instance, increasing $\theta_k^s$ scales the distribution horizontally and stretches the range of the distribution. Similarly, decreasing $\theta_k^s$ scales the distribution vertically and compresses the range of the distribution. If the mean is less than the standard deviation ($\lambda_k^s<1$) the gamma-distribution is exponentially shaped and asymptotic to the vertical and horizontal axes. If the mean is equal to the standard deviation ($\lambda_k^s=1$) the gamma-distribution is exponentially shaped with mean $\theta_k^s$. If the mean is greater than the standard deviation ($\lambda_k^s>1$) then the gamma-distribution is unimodal and skewed-shape, where the skewness decreases in $\lambda_k^s$.

In the next section, the PCC voltage of five PV systems, each with different line impedance, are analyzed and statistical parameters regarding the accuracy of these models are presented.

\section{Real Data and Simulations of Voltage Model}\label{Sectionft}

\begin{figure}[h!]
\centering
\framebox{\parbox{3.2in}{
\includegraphics[width=3in]{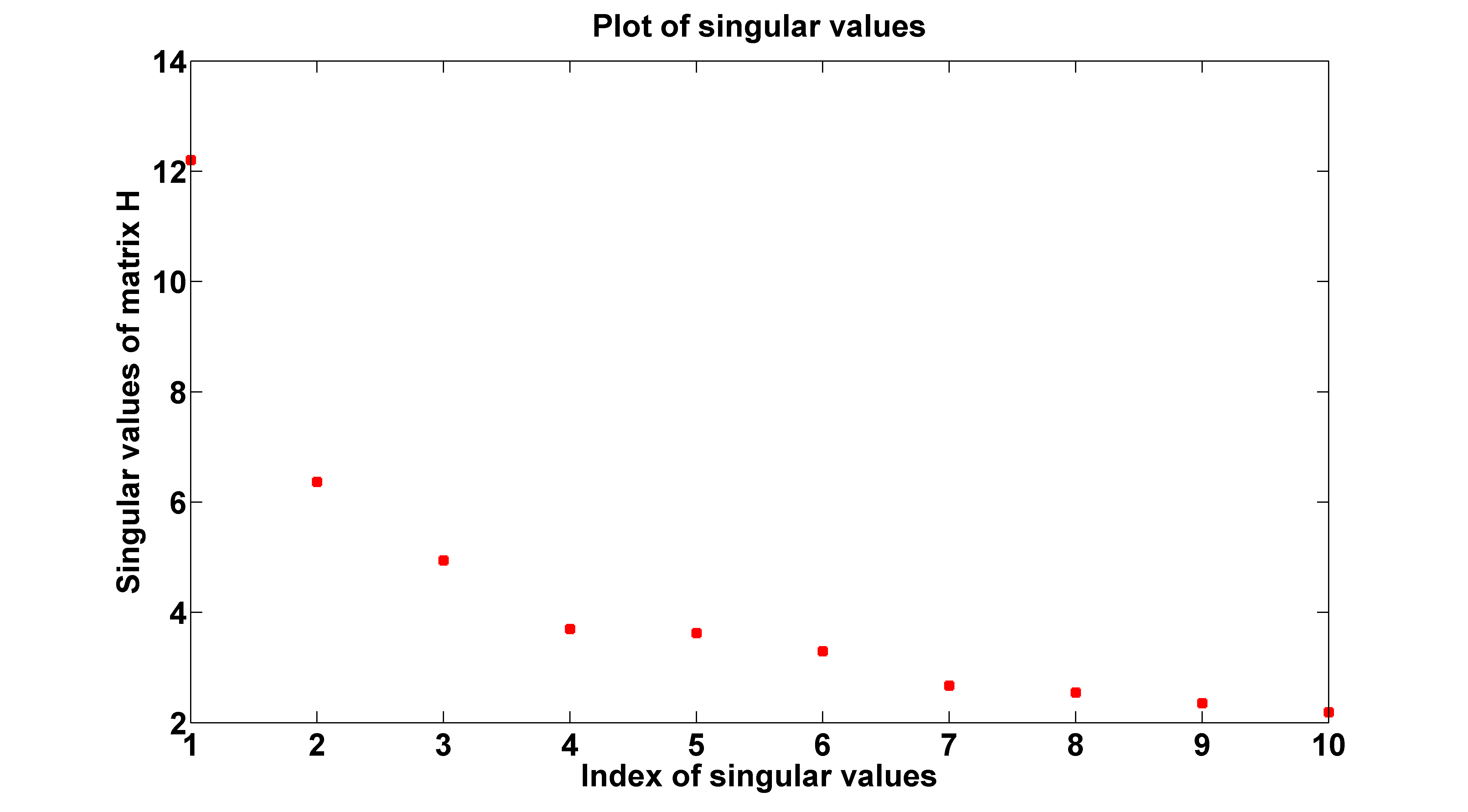}}}
\caption{Singular values of the matrix $H$.}
\label{fig:singularvalues}
\end{figure}

We model the voltage rise caused by PV system output. The analysis is based on the voltage and active power data of PCCs from $10:00$ a.m. to $5:00$ p.m. Let the data sets of the five PV systems be stacked in matrix $H$. The first ten singular values of matrix $H$ are plotted in Figure~\ref{fig:singularvalues}. It is observed that the first three singular values are significantly larger than the remaining singular values.

\begin{figure}[h!]
\centering
\framebox{\parbox{3.2in}{
\includegraphics[width=3in]{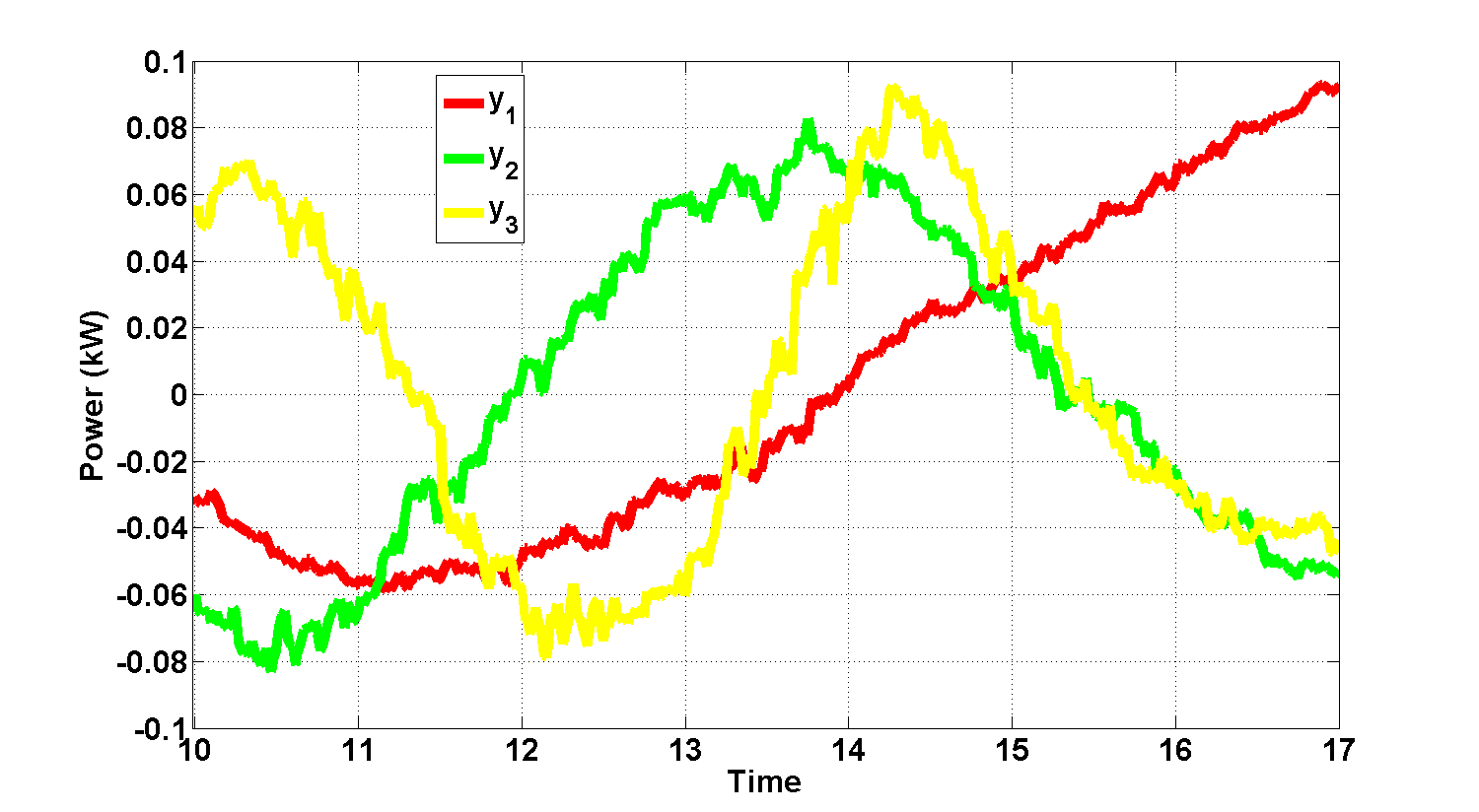}}}
\caption{$y_1$, $y_2$ and $y_3$.}
\label{fig:signals}
\end{figure}

Let $y_1$, $y_2$ and $y_3$ represent the right singular vector in the first, second and third repetitions of the $\ell_1$-SVD Algorithm, as given in Figure~\ref{fig:signals}. $y_1$, $y_2$ and $y_3$ are the normalized bases for active power at the consumer PCC.
It is observed that the basis $y_1$ presents a profile in which consumers inject power into the grid from $10:00$ a.m. to $2:00$ p.m.,
and begin consuming power from the grid after $2:00$ p.m. $y_2$ yields a profile in which consumers inject power into the grid before noon and after $3:30$ p.m. $y_3$ presents a profile with increased fluctuations in the demanded load and PV system output. For the sake of simplicity, $y_3$ is not considered for the grouping of the rows of data matrix $H$. Using $y_1$ and $y_2$, the data matrix $H$ is grouped into three distinct subsets. The subset $\mathcal{C}_1$ ($\mathcal{C}_2$) corresponds to the days that $x_1$ has non-zero entries in the first (second) repetition, on these days the aggregate power is highly correlated with $y_1$ ($y_2$). The subset $\mathcal{C}_3$ corresponds to the days that $x_1$ has zero entries in both the first and second repetitions, and on these days the aggregate power is not correlated with $y_1$ and $y_2$.

The coefficients $\pi_k^s$ and the Gamma distribution parameters $\lambda_k^s$ and $\theta_k^s$ are given in Tables II and IV. From Tables II and IV, it is observed that:

\begin{table}[h!]
\caption{Regression parameter $\beta$}
\label{table:tb0000}
\begin{center}
\begin{tabular}{|c|c|c|c|c|c|}
\hline
  kW &1.9&3.9&7.3&11.6&9.2\\
\hline
   $\beta$  &-0.0113 & -0.0094 & -0.0036 & -0.0018 & -0.0014\\
\hline
\end{tabular}
\end{center}
\end{table}

\begin{table}[h!]\label{tb:tb1000}
\caption{Coefficients $\pi_k^s$}
\begin{center}
\begin{tabular}{|c|c|c|c|c|c|}
\hline
kW & 1.9 & 3.9 & 7.3  & 11.6 & 9.2 \\
\hline
$\pi_1^{+}$   &      0  &       0  &  0.003   &  0.236 &  0.051 \\
\hline
$\pi_1^{-}$   &      0  &       0  &  0.003   &  0.223 &  0.042 \\
\hline
$\pi_2^{+}$   & 0.020  &  0.003  &  0.017   &  0.147 &  0.121 \\
\hline
$\pi_2^{-}$   & 0.014  &  0.003  &  0.018  &  0.138 &   0.110\\
\hline
$\pi_3^{+}$   & 0.528  &  0.558  &  0.496             &  0.136 &  0.348 \\
\hline
$\pi_3^{-}$   & 0.438  &  0.436  &  0.463   &  0.120 &  0.328 \\
\hline
\end{tabular}
\end{center}
\end{table}
\begin{table}[h!]\label{tb:tb2000}
\caption{Gamma distribution parameters $\lambda_k^s,\theta_k^s$ }
\begin{center}
\begin{tabular}{|c|c|c|c|c|c|}
\hline
kW & 1.9 & 3.9 & 7.3 & 11.6 & 9.2 \\
\hline
$u_1^{+}$   &        &         &  1.3,\,\ 0.003  &  1.5,\,\ 0.003  & 1.5,\,\ 0.003 \\
\hline
$u_1^{-}$   &        &         &  1.0,\,\ 0.004   &  1.4,\,\ 0.004  & 1.4,\,\ 0.004 \\
\hline
$u_2^{+}$   & 1.8,\,\ 0.003  &  1.7,\,\ 0.002  &  1.5,\,\ 0.004   &  1.6,\,\ 0.003 & 1.5,\,\ 0.003  \\
\hline
$u_2^{-}$   & 1.1,\,\ 0.007  &  1.5,\,\ 0.002  &  1.4,\,\ 0.004   &  1.4,\,\ 0.004 & 1.3,\,\ 0.003 \\
\hline
$u_3^{+}$   & 1.6,\,\ 0.003  &  1.5,\,\ 0.004  &  1.4,\,\ 0.003  &  1.5,\,\ 0.003 &   1.4,\,\ 0.003\\
\hline
$u_3^{-}$   & 1.1,\,\ 0.006  &  1.1,\,\ 0.007  &  1.3,\,\ 0.004   &  1.4,\,\ 0.003 & 1.4,\,\ 0.003  \\
\hline
\end{tabular}
\end{center}
\end{table}

\begin{itemize}
\item Regression parameter $\beta$ is decreasing in the line impedance.
\item All of the Gamma-distributions are unimodal, because $\lambda_k^s \geq 1$.
\item $\lambda_k^+ \geq \lambda_k^-$ and $\theta_k^+ \leq \theta_k^-$; consequently, $u_k^-$ has a larger range of distribution than $u_k^+$, and $u_k^+$ has a less skew-shaped distribution than $u_k^-$.
\item $\pi_3^+$ coefficients have a higher value than the other coefficients for each of the PV systems with capacity $1.9$, $3.9$, $7.3$ and $9.2$ kW.
$\pi_1^+$ has a higher coefficient value than the other coefficients for the $11.6$ kW PV system.
\end{itemize}

\begin{figure}
\centering
\framebox{\parbox{3in}{
\includegraphics[width=3in]{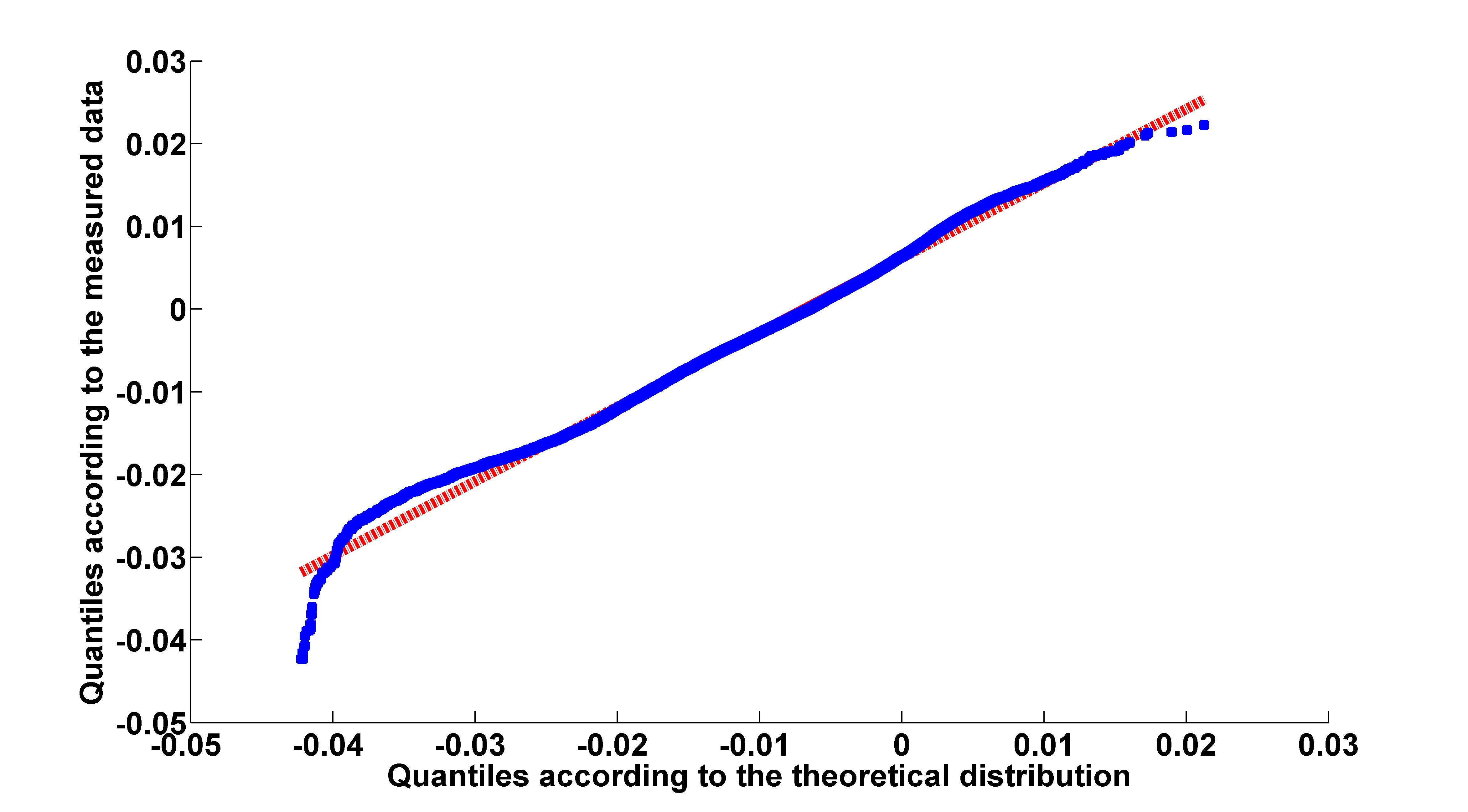}}}
\caption{Q-Q plot with clustering the data. }
\label{fig:qqplot2}
\end{figure}

\begin{figure}
\centering
\framebox{\parbox{3in}{
\includegraphics[width=3in]{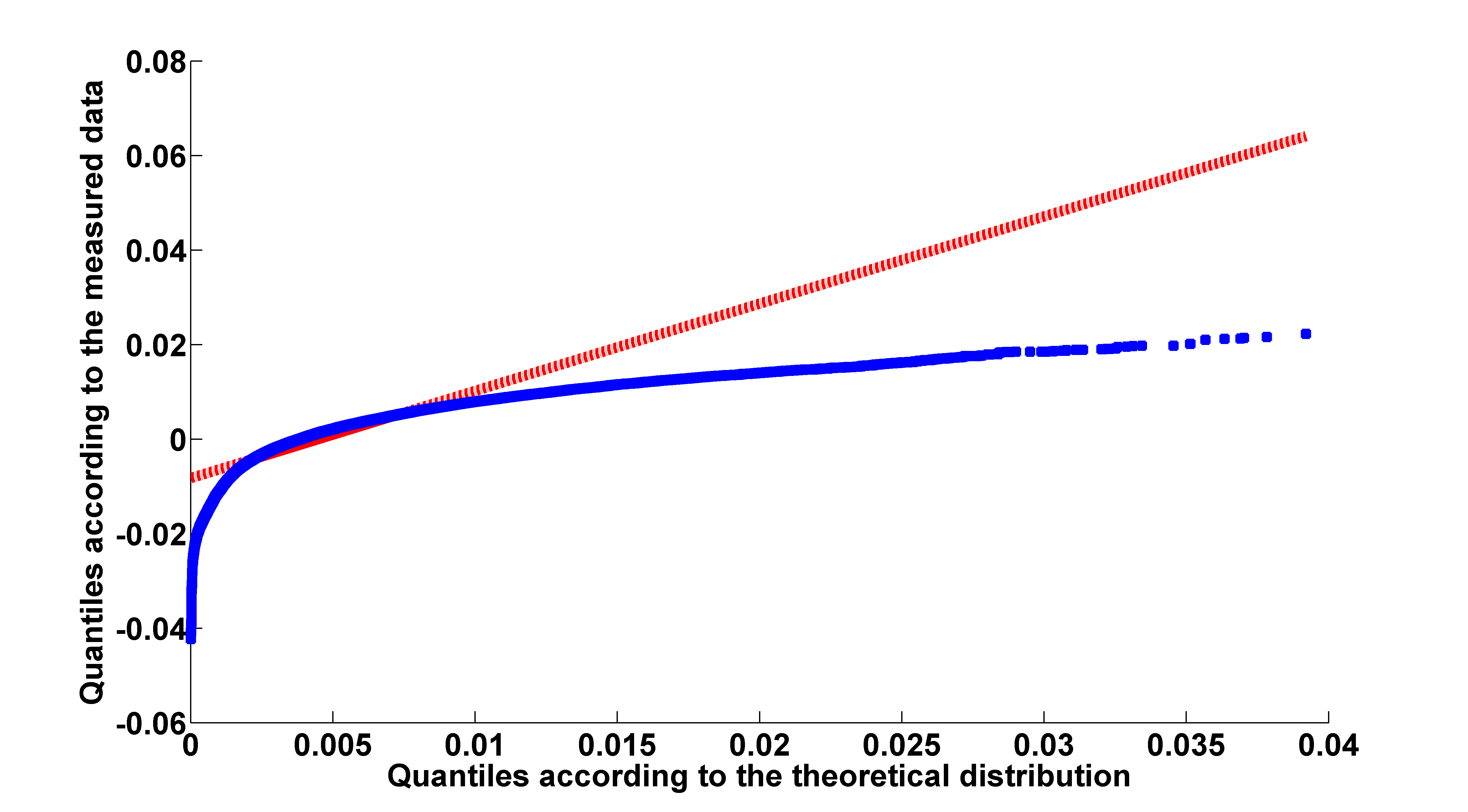}}}
\caption{Q-Q plot without clustering the data. }
\label{fig:qqplot3}
\end{figure}

We use a Quantile-Quantile (Q-Q) plot (\cite{Wilk}) to compare the distribution of the measurement $V-\beta P$ and the theoretical distribution $\sum_{k=1}^{3} \pi_k^+ u_k^+-\pi_k^- u_k^-$. A point $(x,y)$ on the Q-Q plot corresponds to a quantile of the second data set against the same quantile of the first data set. If the sizes of two data sets are equal, then the Q-Q plot is the plot of sorted data of the second data set versus the sorted data of the first data set. If the first and second data sets come from the same distribution then the Q-Q plot will fit close to the identity line $x=y$. Let the first data set be the data generated by the distribution $\sum_{k=1}^{3} \pi_k^+ u_k^+-\pi_k^- u_k^-$ and the second data set be the residue $V-\beta P$. We consider equal sizes when computing both data sets. In Figure~\ref{fig:qqplot2}, the Q-Q plot of the two quantiles set against one another is shown. Now, assume that the first data set is the residue $V-\beta P$ and the second data set is the data generated by a Gamma distribution fitted to data $V-\beta P$, using maximum likelihood estimation and without clustering the data. In Figure~\ref{fig:qqplot3}, the Q-Q plot of the quantiles of the first data set against the quantiles of the second data set is shown. The divergence from the linear line, is evidence that the two data sets come from different distributions. Unlike Figure~\ref{fig:qqplot3}, Figure~\ref{fig:qqplot2} is linear $x=y$, therefore clustering is an effective approach, enabling us to better model the residue of least square regression.

\section{Voltage Regulator}\label{se:vr}
\begin{figure}[h!]
\centering
\framebox{\parbox{3.2in}{
\includegraphics[width=3in]{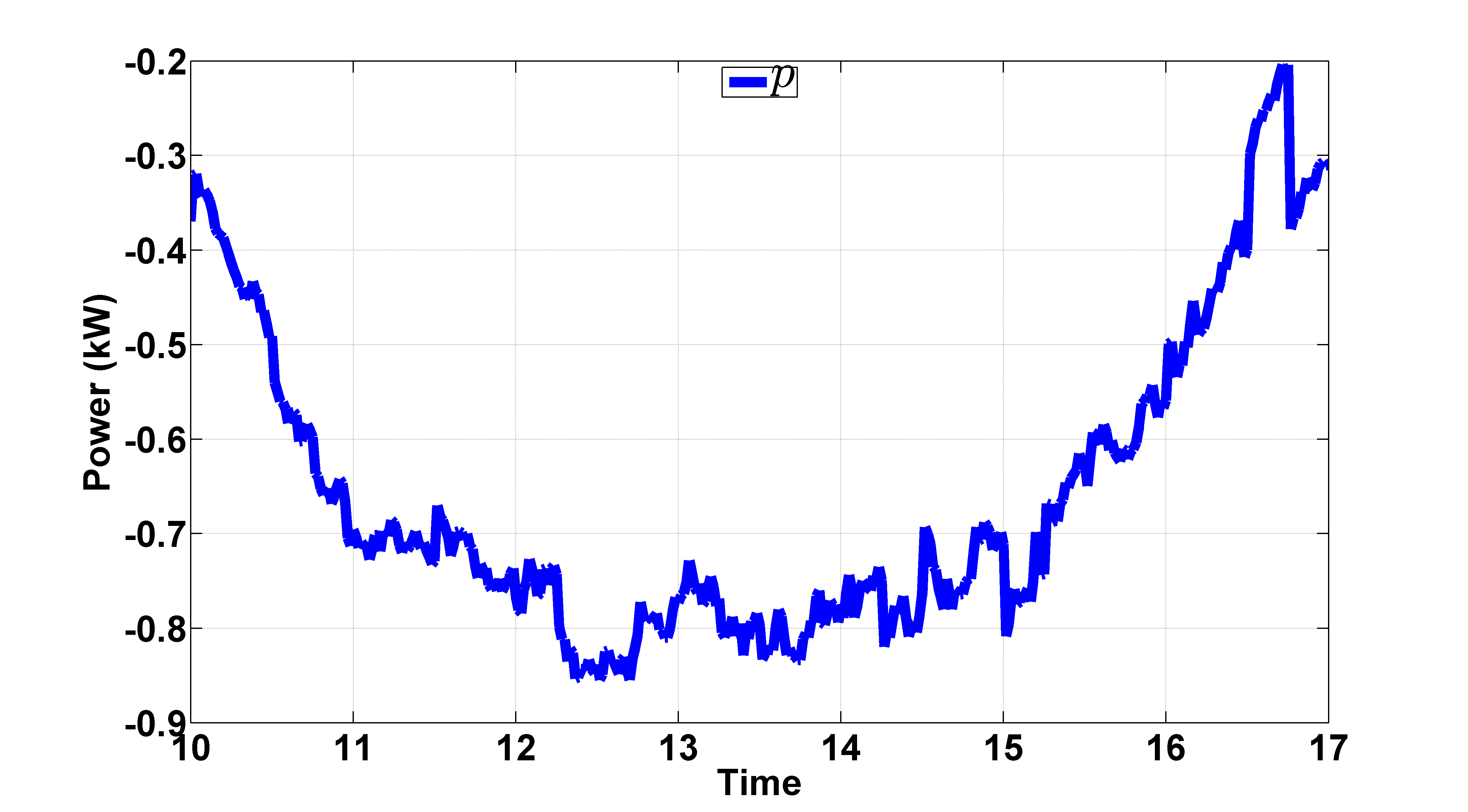}}}
\caption{Average of the measured power $p(t)$.}
\label{fig:avg_power}
\end{figure}

\begin{figure}[h!]
\centering
\framebox{\parbox{3.2in}{
\includegraphics[width=3in]{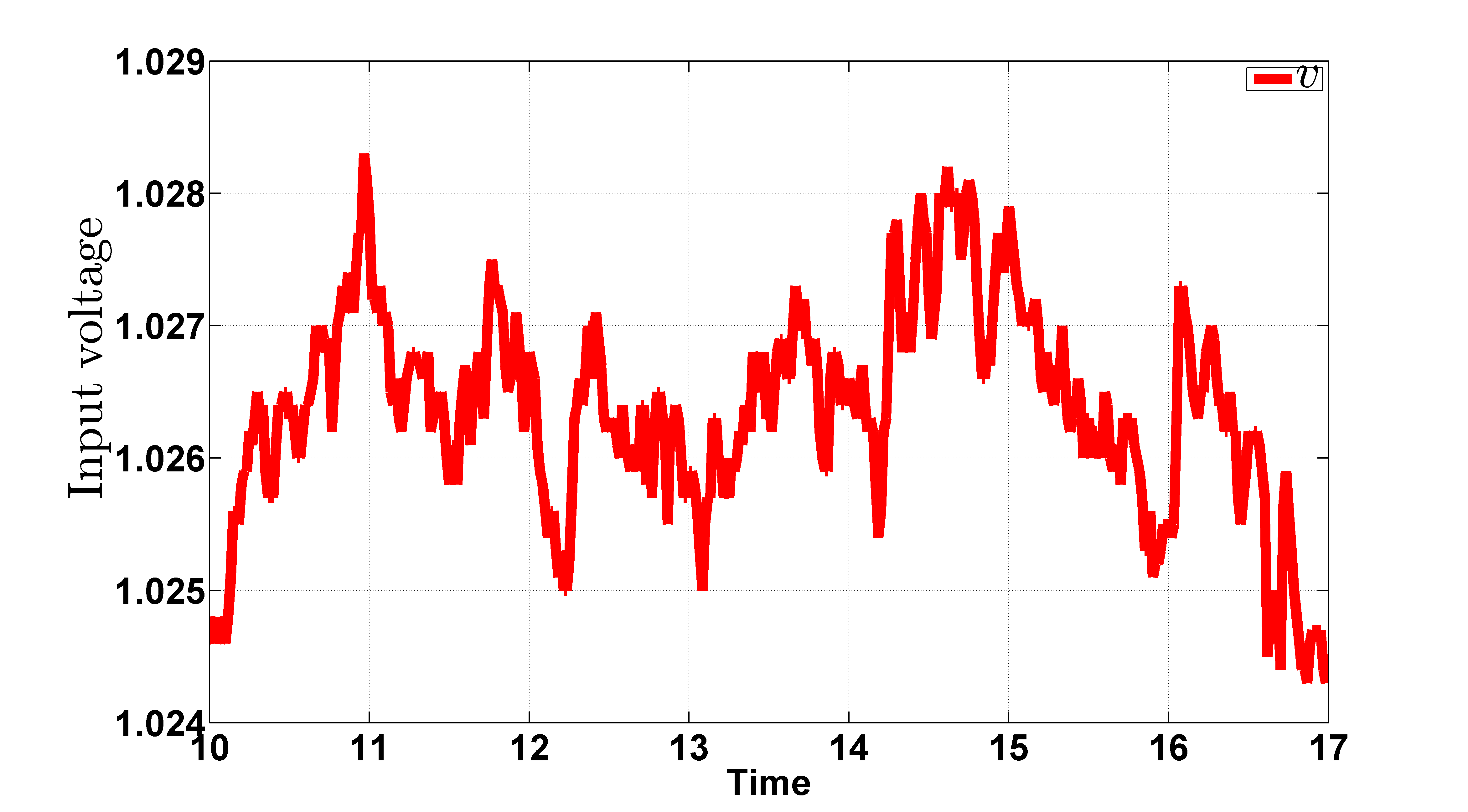}}}
\caption{Average of the measured voltage $v(t)$.}
\label{fig:avg_volt}
\end{figure}

\begin{figure}[h!]
\centering
\framebox{\parbox{3.2in}{
\includegraphics[width=3in]{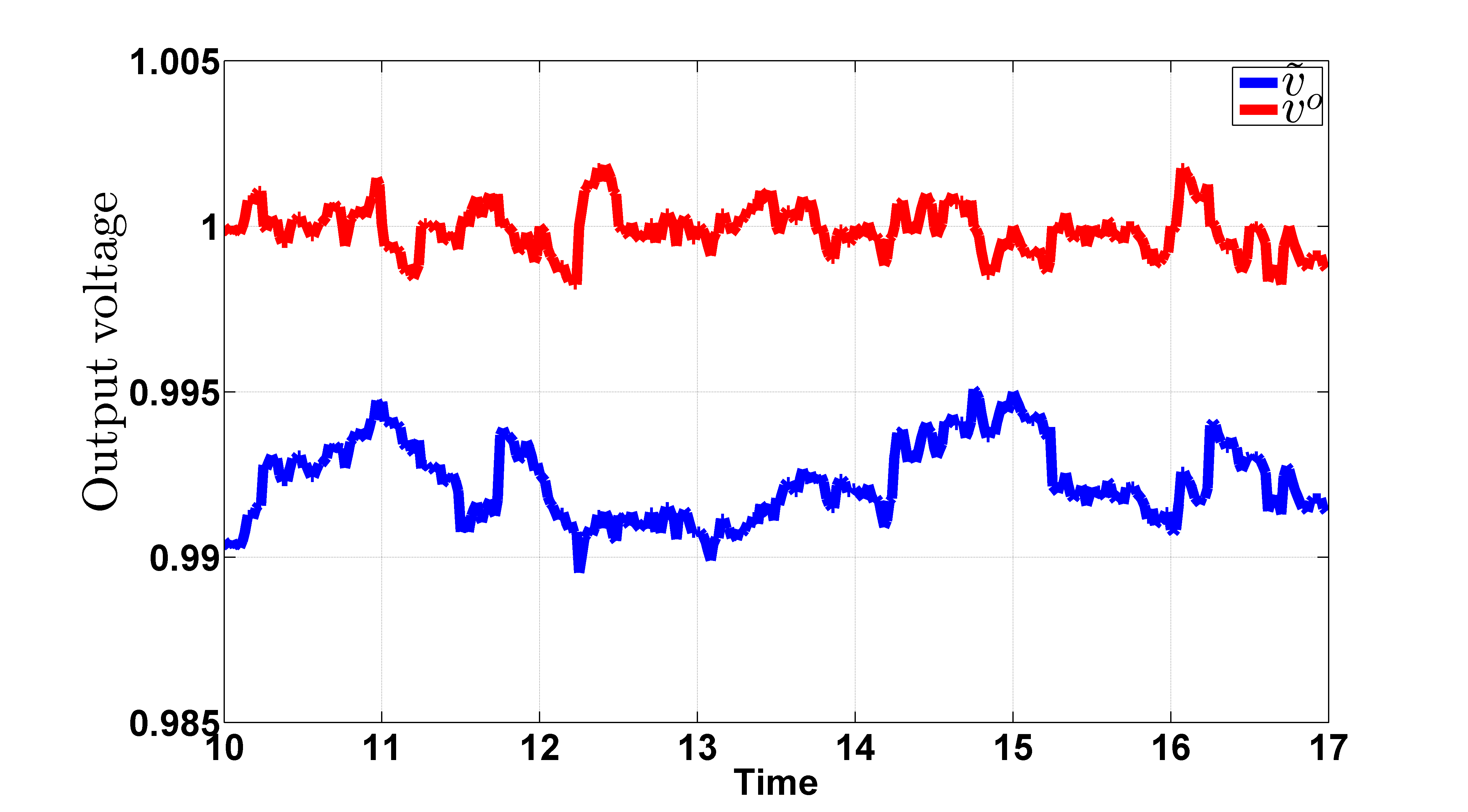}}}
\caption{Output voltage of conventional voltage regulator $v^{o}(t)$ and output voltage of the voltage regulator based on stochastic model $\tilde{v}^{o}(t)$.}
\label{fig:outputs}
\end{figure}

Efficient adjustment of the transformer LTC can alleviate unnecessary inverter tripping and system losses caused by voltage rise. The PV systems on the distribution grids can cause rapid voltage fluctuations. This causes difficulties for the operation of the transformer LTC, which regulates the voltage at a slower timescale. We consider a transformer LTC that regulates the voltage at the PCC at a slow timescale $T$. Let $p(t)$ and $v(t)$ be the average of the measured power and voltage at the PCC from $10:00$ a.m. to $5:00$ p.m., as shown in Figures~\ref{fig:avg_power} and \ref{fig:avg_volt}. The input to the voltage regulator is $v(t)$. The output voltage of the conventional voltage regulator is given as $v^{o}(t)=\frac{1}{v(T)}v(t)$ for all $t$ in $T$, as shown in Figure~\ref{fig:outputs} in red. The LTC value for the conventional voltage regulator is defined as $LTC(t)=\frac{1}{v(T)}$, for all $t$ in $T$. We develop a new voltage regulator that operates based on measured voltage, power, and the stochastic model developed in the previous section. Let $n=\sum_{k=1}^{3} \pi_k^+ u_k^+-\pi_k^- u_k^-$. $F_n$ and $\acute{F}_n$ are respectively the Cumulative Distribution Function (CDF) of $n$ and the derivative of the CDF. Let $E$ be the expected value over $n$ and $\Delta$ be a constant. Let $v_d(T)=v(T)-\beta p(T)$ and $n_1$ be defined as
\begin{align}
n_1=\min\{z|  F_n(z) \geq & F_n(v_d(T)) -\Delta \acute{F}_n\big(v_d(T)\big)\}.
\end{align}
The value of voltage $v_d(T)$ will not exceed $n_1$ with probability $F_n(v_d(T)) -\Delta \acute{F}_n\big(v_d(T)\big)$. Let $n_2$ be defined as
\begin{align}
n_2=\max\{z|  F_n(z) < & F_n(v_d(T)) +\Delta \acute{F}_n\big(v_d(T)\big)\}.
\end{align}
The value of voltage $v_d(T)$ will exceed $n_2$ with probability $1-F_n(v_d(T)) -\Delta \acute{F}_n\big(v_d(T)\big)$. Let $E$ denote the expected value over $n$. We define $\gamma_t$ as the conditional expectation of $n$, given $n$ is greater than $n_1$ and less than $n_2$
\begin{align}
\gamma_t=E \Big[n \Big{|} \,\  n_1 < n \leq n_2 \Big]
\end{align}
for all $t$ in $T$. The output voltage of the voltage regulator is given as $\tilde{v}^{o}(t)=\frac{1}{\beta p(T)+\gamma_t}v(t)$, as shown in Figure~\ref{fig:outputs} in blue. The LTC value for this voltage regulator is defined as $LTC(t)=\frac{1}{\beta p(T)+\gamma_t}$, for all $t$ in $T$.

\begin{figure}[h!]
\centering
\framebox{\parbox{3.2in}{
\includegraphics[width=3in]{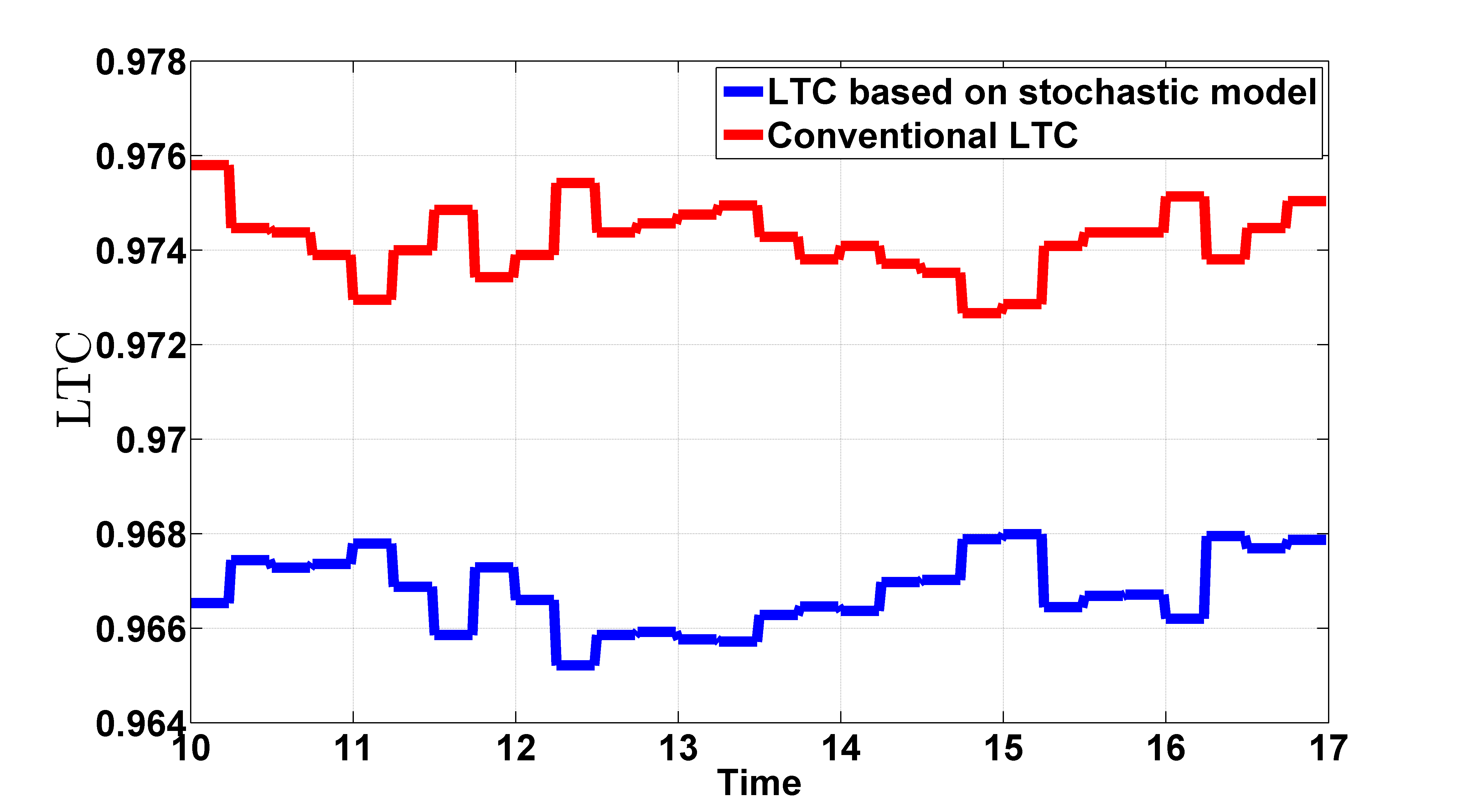}}}
\caption{LTC position.}
\label{fig:LTC}
\end{figure}

It is observable from Figure~\ref{fig:outputs}, that the voltage regulator incorporating the stochastic model, lowers the voltage more than the conventional voltage regulator at times of high PV output. Consequently, the developed LTC controller decreases the power consumption of consumers and increases the potential capacity for PV systems on distribution grids without voltage violations. The LTC switches are shown in the Figure~\ref{fig:LTC}. It is observed that the LTC that operates based on the stochastic model has a lower value of LTC variation $(\sum_{t} |LTC(t)-LTC(t-1)|)$, $0.0149$ versus $0.0177$. Therefore, the stochastic model can also contribute to improving the transformer lifetime.

\section{CONCLUSIONS} \label{SectionVII}

In this work, the effect of injected power on voltage rise is quantified. The voltage at each PCC is affected by neighboring loads and PV system outputs. The voltage variation caused by injected power of solar PV systems, is modeled as a linear combination of Gamma random variables. The $\ell_1$-SVD is used to find the bases of the active power data of the PCC of consumers on a radial distribution grid.
Each subset of the data is modeled as a Gamma random variable through maximum likelihood estimation. The goal of clustering the data is to model the sparse voltage rise in the voltage data. Since the data sets consist of different features, fitting a model to the data can be complicated in the absence of clustering. The algorithm developed in this work provides a tool to model the voltage rises caused by PV systems, taking into account the sparse events in the voltage rises. Least squares regression is not able to model high voltage rises that occur over short time intervals. Clustering data sets, especially for big data sets, enables us to model voltage rises that are sparse in the data sets. The stochastic model is used for developing a voltage regulator that lowers the voltage more effectively than conventional voltage regulators at times of high PV generation. Therefore, this model increases the potential capacity of PV solar systems on distribution grids. This model also decreases the value of LTC. Therefore it helps to increase the lifetime of the voltage regulator.

\addtolength{\textheight}{-12cm}   % This command serves to balance the column lengths

%%%%%%%%%%%%%%%%%%%%%%%%%%%%%%%%%%%%%%%%%%%%%%%%%%%%%%%%%%%%%%%%%%%%%%%%%%%%%%%%

%%%%%%%%%%%%%%%%%%%%%%%%%%%%%%%%%%%%%%%%%%%%%%%%%%%%%%%%%%%%%%%%%%%%%%%%%%%%%%%%

%%%%%%%%%%%%%%%%%%%%%%%%%%%%%%%%%%%%%%%%%%%%%%%%%%%%%%%%%%%%%%%%%%%%%%%%%%%%%%%%

%\section*{ACKNOWLEDGMENT}

%This project is funded by US Department of Energy grant number XXXXXX and National Science Foundation grant number XXXXXX.


\begin{thebibliography}{99}

\bibitem{conferenceversion}
A. Zeinalzadeh, R. Ghorbani, J. Yee, Stochastic model of voltage variations in the presence of photovoltaic systems, IEEE American Control Conference (ACC), Boston, MA, USA, July 6-8, 2016.
\bibitem{c15} T. Stetz, F. Marten, and M. Braun, Improved low voltage grid integration of photovoltaic systems in Germany, IEEE Transaction on  Sustainability Energy, vol. 4, no. 2, pp. 534-542, 2013.
\bibitem{c17} E. Demirok, D. Sera, R. Teodorescu, P. Rodriguez, and U. Borup, Evaluation of the voltage support strategies for the low voltage grid connected pv generators, in Proc. IEEE Energy Conversion Congress and Exposition (ECCE), pp. 710-717, 2010.
\bibitem{c177} A. Canova, L. Giaccone, F. Spertino, and M. Tartaglia, Electrical impact of photovoltaic plant in distributed network, IEEE Transactions on Industry Applications, vol. 45, no. 1, pp. 341-347, 2009.
\bibitem{c1777} M. Thomson and D. G. Infield, Network power-flow analysis for a high penetration of distributed generation, IEEE Transactions on Power Systems, vol. 22, no. 3, pp. 1157-1162, 2007.
\bibitem{c17777} R. A. Walling, R. Saint, R. C. Dugan, J. Burke, and L. A. Kojovic, Summary of distributed resources impact on power delivery systems, IEEE Transactions on Power Delivery, vol. 23, no. 3, pp. 1636-1644, 2008.

\bibitem{c1} R. J. Broderick, J. E. Quiroz, M. J. Reno, A. Ellis, J. Smith, and R. Dugan, Time series power flow analysis for distribution connected pv generation, Sandia National Laboratories SAND-0537, 2013.
\bibitem{c10} T. Stetz, W. Yan, and M. Braun, Voltage control in distribution systems with high level PV-penetration, in 25th
    European PV Solar Energy Conference, Valencia, Spain, 2010.
\bibitem{c7} T. Verschueren, K. Mets, B. Meersman, M. Strobbe, C. Develder, L. Vandevelde, Assessment and mitigation of voltage violations by solar panels in a residential distribution grid, IEEE International Conference on Smart Grid Communications (SmartGridComm), pp. 540-545, 2011.
\bibitem{c77}
A. Zeinalzadeh and V. Gupta, Minimizing risk of load shedding and renewable energy curtailment in a microgrid with energy storage, ArXiv
e-prints, arXiv:1611.08000, Nov. 2016.
\bibitem{c8} B. Perera, P. Ciufo, and S. Perera., Point of common coupling (PCC) voltage control of a grid-connected solar photovoltaic (PV) system, Annual Conference of the IEEE Industrial Electronics Society. Vienna, Austria. Jan. 2013.
\bibitem{Alam}
M. J. E. Alam , K. Muttaqi and D. Sutanto, Distributed energy storage for mitigation of voltage-rise impact caused by rooftop solar PV, IEEE Transactions on Power Systems, vol. 28, no. 4, pp. 3874-3884, 2013.
\bibitem{Ueda}
Y. Ueda, K. Kurokawa, T. Tanabe, K. Kitamura, K. Akanuma, M. Yokota, and H. Sugihara, Study on the over voltage problem and battery
operation for grid-connected residential PV systems, 22nd European Photovoltaic Solar Energy Conference, 3-7 September 2007.
\bibitem{Sera}
E. Demirok, D. Sera, R. Teodorescu, P. Rodriguez, and U. Borup, Clustered PV inverters in LV networks: An overview of impacts and
comparison of voltage control strategies, in Electrical Power and Energy Conference (EPEC),  IEEE, pp. 1-6, 2009.
\bibitem{ashkanreactive}
A. Zeinalzadeh, R. Ghorbani, E. Reihani, Optimal power flow problem with energy storage voltage and reactive power control, The 45th ISCIE International Symposium on Stochastic Systems Theory and Its Applications, Okinawa, 2013.
\bibitem{Ali}
S. Ali, N. Pearsall and G. Putru, Impact of high penetration level of grid-connected photovoltaic systems on the uk low voltage distribution network, International Conference on Renewable Energies and Power Quality, Santiago de Compostela, pp. 1-4, 2012.
\bibitem{Wilk}
M.B. Wilk, R. Gnanadesikan, Probability plotting methods for the analysis of data, Biometrika, vol. 55, no. 1, pp. 1-17, 1968.

\end{thebibliography}
\end{document}